\documentclass{article}
\usepackage{amsmath}[1996/11/01]
\usepackage{amssymb,amsthm,amsxtra}
\usepackage{epic,eepic}
\usepackage{epsfig}

\def\[#1\]{\begin{equation}#1\end{equation}}
\makeatletter
\def\beq{%
   \relax\ifmmode
      \@badmath
   \else
      \ifvmode
         \nointerlineskip
         \makebox[.6\linewidth]%
      \fi
      $$
   \fi
}
\def\eeq{%
   \relax\ifmmode
      \ifinner
         \@badmath
      \else
         $$
      \fi
   \else
      \@badmath
   \fi
   \ignorespaces
}

\def\enddisplaymath{\eeq\global\@ignoretrue}
\makeatother

\newtheorem{thm}{Theorem}

\newtheorem{lem}[thm]{Lemma}
\newtheorem{prop}[thm]{Proposition}
\theoremstyle{remark}
\newtheorem*{rem}{Remark}

\theoremstyle{definition}

\numberwithin{equation}{section}
\numberwithin{thm}{section}

\DeclareMathOperator{\Exp}{\mathbb E}

\DeclareMathOperator{\Prob}{\mathbb P}

\DeclareMathOperator{\pathcon}{P}

\newcommand{\C}{\mathbb C}
\newcommand{\R}{\mathbb R}

\begin{document}

\title{{\bf Random vicious walks and random matrices}}
\author{{\bf Jinho Baik}\footnote{
Princeton University and Institite for Advanced Study,
New Jersey, jbaik@math.princeton.edu}}

\date{December 29, 1999}
\maketitle

\begin{abstract}
Lock step walker model is a one-dimensional integer lattice walker model 
in discrete time.
Suppose that initially there are infinitely many walkers 
on the non-negative even integer sites. 
At each tick of time, each walker moves either to its left or to its right 
with equal probability.
The only constraint is that no two walkers can occupy the same site at the 
same time.
It is proved that in the large time limit, a certain conditional probability 
of the displacement of the leftmost walker is identical to the limiting 
distribution of the properly scaled largest eigenvalue of a random GOE matrix 
(GOE Tracy-Widom distribution).
The proof is based on the bijection between path configurations and 
semistandard Young tableaux established recently 
by Guttmann, Owczarek and Viennot.
Statistics of semistandard Young tableaux is analyzed using 
the Hankel determinant expression for the probability from the work of 
Rains and the author. 
The asymptotics of the Hankel determinant
is obtained by 
applying the Deift-Zhou 
steepest-descent method to 
the Riemann-Hilbert problem for the related orthogonal polynomials.

\end{abstract}


\section{Introduction}\label{sec-intro}

In \cite{Fisher}, two types of random vicious 
walkers models, \emph{random turn walker model} and 
\emph{lock step walker model}, are considered.
In these models, walkers are on one-dimensional integer lattice, 
and time is discrete. 
For their applications and earlier results, see, for example,
\cite{Fisher, For1, For2, For3, GOV, BrEsOw1, BrOw1, 
forrester99} and references therein.
In this paper, we present results on lock step model showing a relation 
to random matrix theory.
For random turn walker model, see \cite{forrester99, BR3} and 
discussions following Theorem \ref{mainthm} below.

\begin{figure}[ht]
 \centerline{\epsfig{file=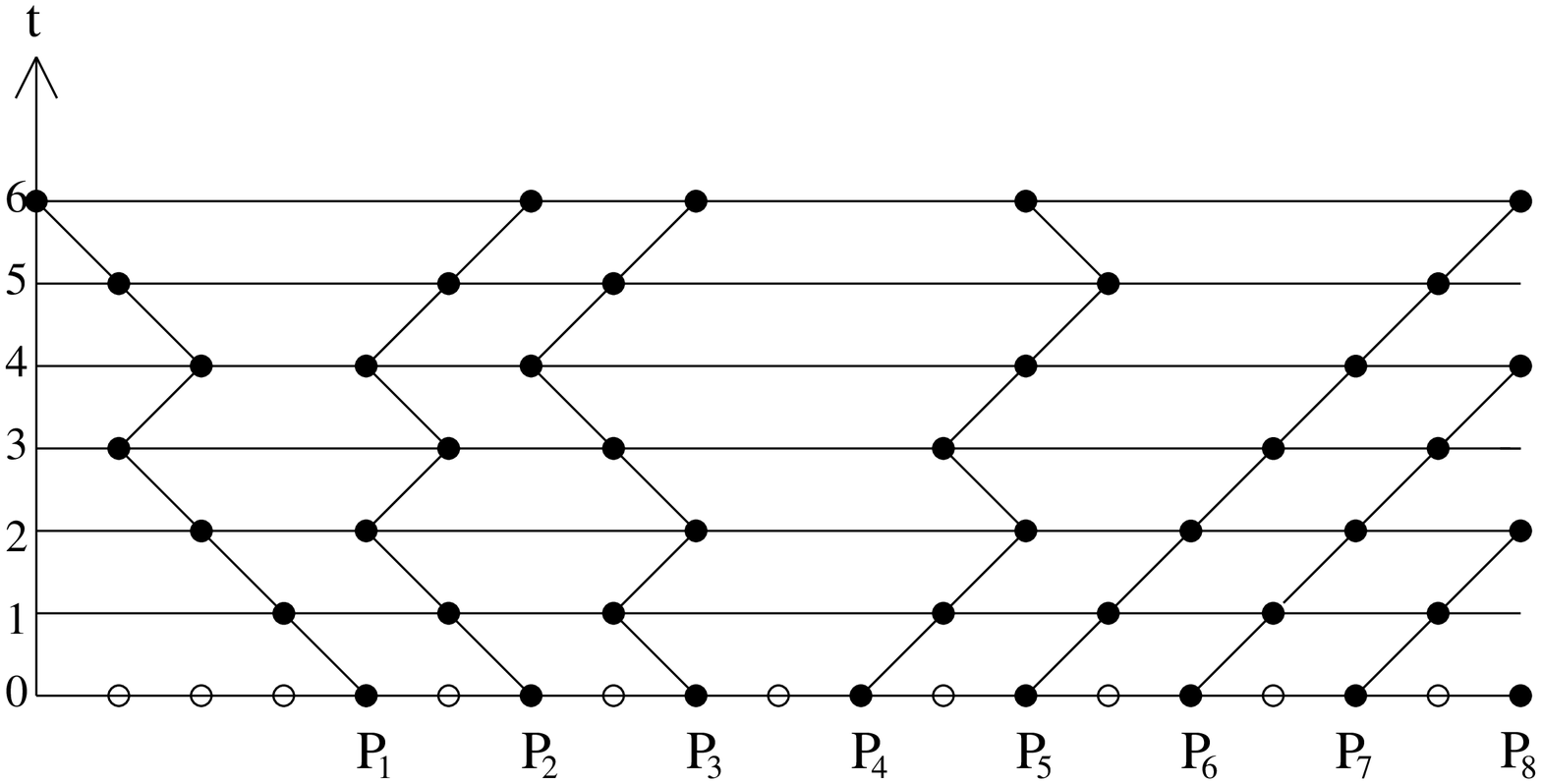, width=11cm}}
 \caption{vicious lock step walkers}
\label{fig-walk}
\end{figure}

At time $t=0$, infinitely many walkers are located at 
the sites $\{0,2,4,6,\cdots\}$.
We label the walkers by $P_1,P_2,P_3,\cdots$ from the left to the right.
In the lock step model, at each time $t=n$, 
all the particles move either to their right 
or to their left with equal probability.
The only constraint is that no two particles can occupy the same site 
at the same time.
This is why the walkers are called ``vicious''.
One typical path configuration is shown in Figure \ref{fig-walk}.

\begin{figure}[ht]
 \centerline{\epsfig{file=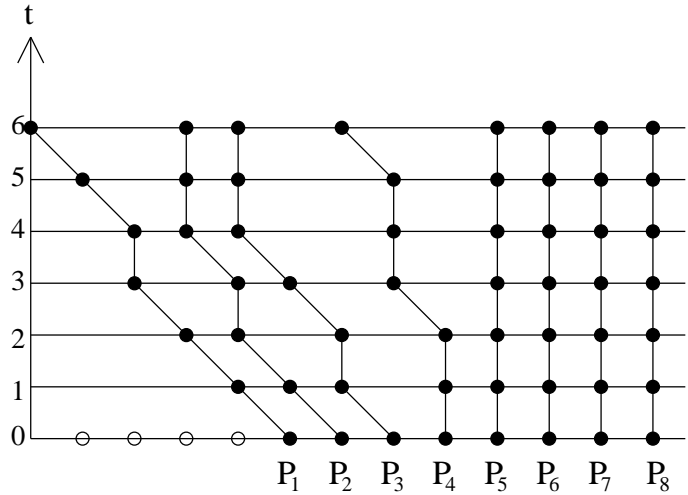, width=9cm}}
 \caption{exclusion process}
\label{fig-newwalk}
\end{figure}

This model can also be thought of as a certain 
totally asymmetric exclusion process in discrete time.
Initially there are infinitely many particles at $\{1,2,3,\cdots\}$. 
A particle is called \emph{left-movable} if its left site is vacant.
Particles $P_{j+1}, P_{j+2}, \cdots, P_{k}$ are called 
\emph{successors} of a particle $P_j$ at a certain time 
if they are next to each other 
in the order of the indices.
At each (discrete) time step, a left-movable particle 
either moves to its left site \emph{together} with 
arbitrarily taken number of its successors, 
or stays at the same site with equal probability.
It is easy to see that this process is equivalent to the above 
lock step model ; right move of lock step corresponds staying 
at the same site in the exclusion process. 
Figure \ref{fig-newwalk} represent an exclusion process 
equivalent to the lock step model in Figure \ref{fig-walk}.

Suppose that during $N$ time steps, total $k$ left moves are 
made by all the particles.
In the example of Figure \ref{fig-walk},  
$N=6$ and $k=14$.
We denote by $\pathcon(N,k)$ the set of all path configurations 
during $N$ time steps with total $k$ left moves.
Then each configuration has equal probability, $1$ over the cardinal 
of $\pathcon(N,k)$.
Hence our probability space is $\pathcon(N,k)$ with uniform probability
given by $1/|\pathcon(N,k)|$.
We denote by $L_j(N,k)(\pi)$ the number of left moves made by 
the particle $P_j$ in a path $\pi\in\pathcon(N,k)$.
We are interested in the limiting statistics of the random variables 
$L_j(N,k)$ as $N,k\to\infty$.

To state the main result, we need a definition.
Let $u(x)$ be the solution of the differential equation 
\begin{equation}\label{as10}
  u_{xx}=2u^3+xu, \qquad 
  u(x)\sim -Ai(x)\quad\text{as}\quad x\to +\infty,
\end{equation}
where $Ai$ is the Airy function.
The above equation is called Painlev\'e II equation. 
It is known that there is unique global solution to \eqref{as10} 
(see, e.g. \cite{BDJ} and references in it).
Define the function, called the 
\emph{GOE Tracy-Widom distribution function}, by  
\begin{equation}
  F_1(x)= \exp\big\{-\frac12\int_x^\infty(s-x)(u(s))^2ds
+\frac12\int_x^\infty u(s)ds\big\}.
\end{equation}
This is indeed a distribution function, and 
the decay rate is given by 
\begin{eqnarray}
  F_1(x) &=& 1+O(e^{-cx^{3/2}}), \qquad \text{as $x\to +\infty$,}\\
  F_1(x) &=& O(e^{-c|x|^{3}}), \qquad\qquad \text{as $x\to -\infty$,}
\end{eqnarray}
for some $c>0$ (see, e.g. (2.11)-(2.14) of \cite{BR2}).
In \cite{TW2}, Tracy and Widom 
proved that $F_1$ is the limiting distribution of the 
properly centered and scaled largest eigenvalue of a random matrix 
taken from a Gaussian orthogonal ensemble.
The subscription $1$ in $F_1$ is a general convention : 
there are also GUE and GSE Tracy-Widom distribution 
functions $F_2$ and $F_4$ \cite{TW1, TW2}.
One can find general discussion for random matrices in 
\cite{Mehta, deift}.

Now the main theorem is 
\begin{thm}\label{mainthm}
For fixed $0<t<1$, let
\begin{equation}\label{eq;walk1}
  \eta(t)=\frac{2t}{1+t}, \qquad
\rho(t)=\frac{\bigl(t(1-t)\bigr)^{1/3}}{1+t}.
\end{equation}
Let $F_1(x)$ be the GOE Tracy-Widom distribution function.
Under the condition that in $N$ time steps total $k$ left moves
are made, the (conditional) probability distribution of the number $L_1(N,k)$
of left moves made by the leftmost particle satisfies
\begin{equation}\label{eq;walk2}
  \lim_{N\to\infty}
\Prob\biggl(\frac{L_1(N,k)-\eta(t)N}
{\rho(t)N^{1/3}}\le x\biggr)
= F_1(x), \qquad 
\text{when $\ \ k=\frac{t^2}{1-t^2}N^2+o(N^{4/3})$}.
\end{equation}
Also all the moments of the scaled random variable converge 
to the corresponding moments of the limiting random variable.
\end{thm}

In other words,
in the large $N$ limit, the (conditional) fluctuation of the displacement 
of the leftmost particle in lock step model is identical to the 
fluctuation of the largest eigenvalue of a random GOE matrix.
Naturally we expect that the $k^{\text{th}}$ particle corresponds 
to the $k^{\text{th}}$ eigenvalue of random GOE matrix.

\begin{figure}[ht]
 \centerline{\epsfig{file=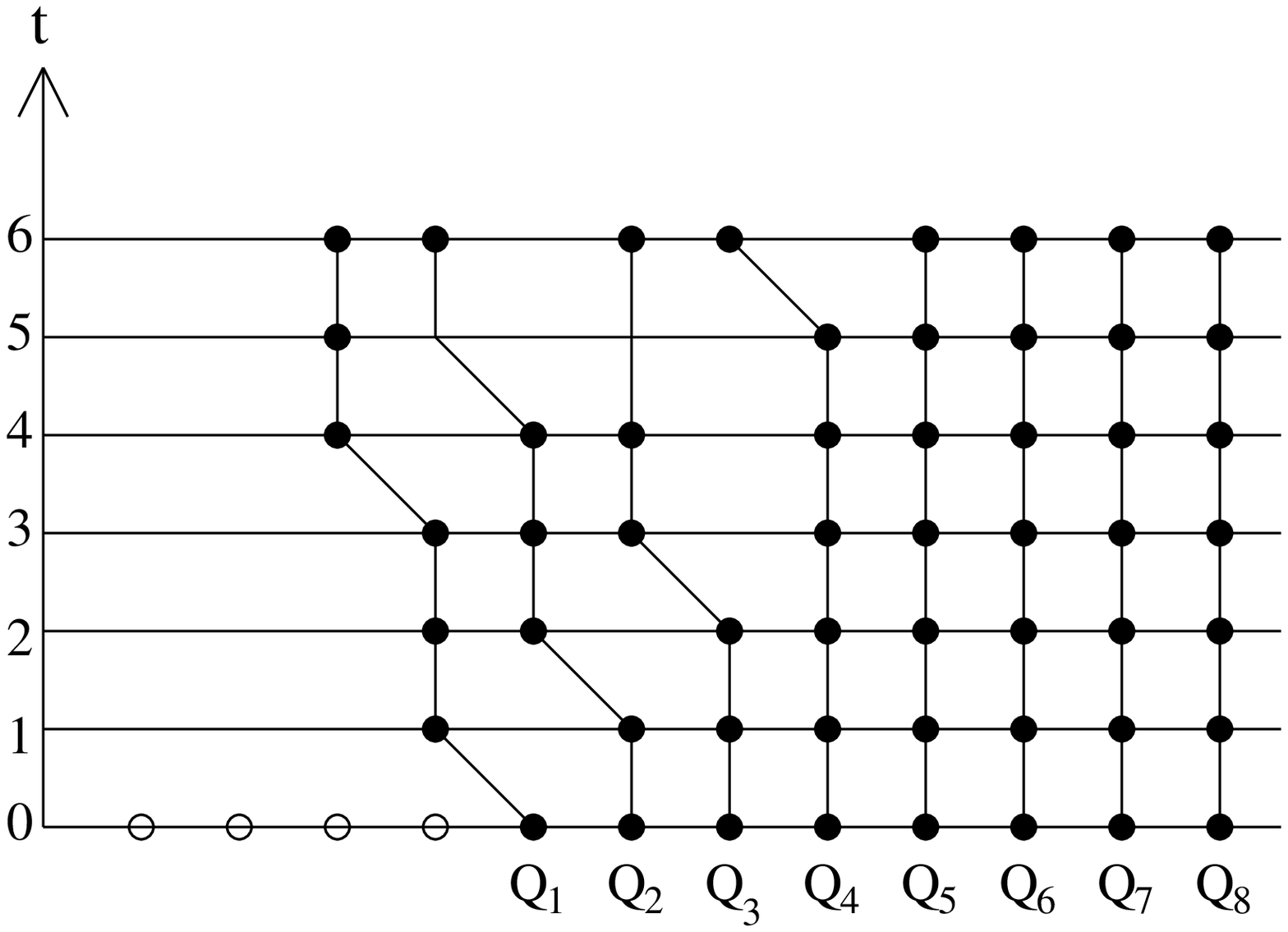, width=9cm}}
 \caption{random turn walker model}
\label{fig-syt}
\end{figure}

It is interesting to compare the above result with 
the results for random turn walker model.
Initially there are infinitely many walkers $Q_1,Q_2,Q_3,\cdots$ 
at the position $\{1,2,3,\cdots\}$.
We again call a walker \emph{left-movable} if its left site 
is vacant.
At each time, we select \emph{one} walker at random among 
left-movable walkers, and move it to its left site.
Hence there is one and only one movement at each time and 
all the movements are to the left.
An example of random turn walker path configuration is in 
Figure \ref{fig-syt}.
Let $X_j(N)$ be the displacement of the $j^{\text{th}}$ walker after 
$N$ time step.
It is shown in \cite{BR3} that for $j=1,2$, we have 
\begin{equation}\label{eq-newin1}
  \lim_{N\to\infty} \Prob\biggl(\frac{X_j(N)-2\sqrt{N}}{N^{1/6}}\le x \biggr)
= F_1^{(j)}(x), \qquad j=1,2,
\end{equation}
where $F_1^{(j)}$ is the limiting distribution of the (scaled) $j^{\text{th}}$ 
largest eigenvalue of a random GOE matrix.
Especially we have $F_1(x)=F_1^{(1)}(x)$.
On the contrary, if we assume that the walkers move 
to their left in the first $N$ time steps, and then move to 
their right in the next $N$ time steps so that at the end 
walkers come back to their original positions, then 
we obtain the GUE Tracy-Widom distribution in the limit \cite{forrester99}. 
Indeed in this case, a lot more are known.
The general $j^{\text{th}}$ row statistics and also 
the correlation functions
converge to the corresponding quantities of random GUE matrix 
in the limit \cite{forrester99}.

\begin{figure}[ht]
 \centerline{\epsfig{file=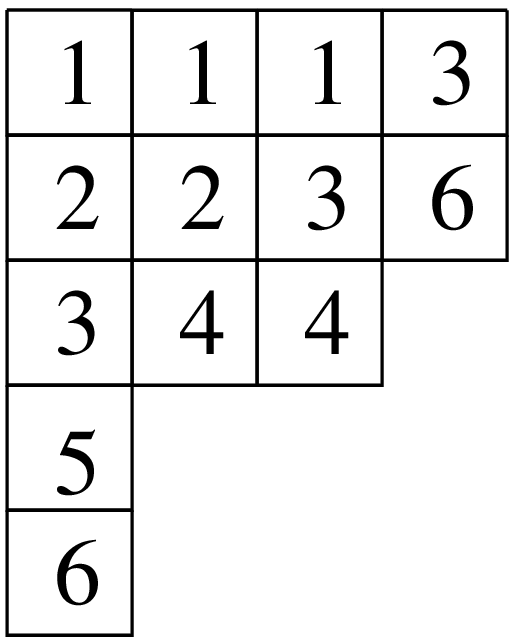, width=3.4cm}}
 \caption{semistandard Young tableau}
\label{fig-tabl}
\end{figure}

\bigskip
The first step to prove the above theorem is to map the path statistics 
into tableaux statistics following \cite{GOV}.
By definition, a semistandard Young tableau (SSYT) 
is an array of positive integers 
top and left adjusted as in Figure \ref{fig-tabl} so that the numbers 
in each row increase weakly and the numbers in each column increase strictly.
A reference for tableaux is \cite{Stl}, and we freely use the notations in it.
In \cite{GOV}, Guttmann, Owczarek and Viennot established a simple bijection 
between path configurations of lock step model and the set of SSYT : 
for a path configuration, we write down the time steps 
at which the $j^{\text{th}}$ 
particle made movement to its left on the $j^{\text{th}}$ column. 
Hence the top row is the array of time steps the particles made first 
movement to their left, the second row is the array of time steps 
the particles made their second movement to their left, and so on.
If we draw boxes around each number, the result is 
a SSYT.
See figure \ref{fig-tabl} for the tableau corresponding to the path 
configuration of Figure \ref{fig-walk}.
This map is a bijection between $\pathcon(N,k)$ and 
the set of SSYT of size $k$ with fillings  
taken from $\{1,2,\cdots,N\}$.
Moreover, under this bijection, $L_j(N,k)$ is equal to 
the number of boxes in the $j^{\text{th}}$ column of the corresponding 
SSYT.
Therefore the statistics of $L_j$ is identical to the 
$j^{\text{th}}$ column statistics of random tableaux.

After the work of Guttmann, Owczarek and Viennot,
Forrester in \cite{forrester99} observed that similar bijection 
can be established between path configuration of random turn model 
and the set of standard Young tableaux (SYT).
The above result \eqref{eq-newin1} is obtained based on this 
bijection and the recent work \cite{BR2} on statistics of 
the first row of random SYT.
Also as mentioned above, if we assume that the walkers move
to their left in the first $N$ time steps, and then move to
their right in the next $N$ time steps so that at the end
walkers come back to their original positions,
the limiting fluctuation is not GOE type, but GUE type.
This difference comes from the fact that in this case, 
the path configuration is in bijection with the \emph{pairs} of SYT.
The statistics of pairs of SYT is more well studied than that of single SYT
\cite{BDJ, BDJ2, Ok, BOO, kurtj:disc}, and we have stronger results 
mentioned earlier.

\bigskip
This paper is organized as follows.
In Section \ref{sec-walk}, using the result of \cite{BR1}, 
we express the generating 
function for the probability of the first column of random 
tableaux in terms of a Hankel determinant. 
It is a general fact that Hankel determinant is 
related to orthogonal polynomials on the unit circle. 
The asymptotics of orthogonal polynomials is obtained via Riemann-Hilbert 
problem and summarized separately in Section \ref{sec-rhp}.
We can obtain the limiting statistics of the first column 
from the knowledge of Hankel determinant asymptotics. 
This work occupies the second half of Section \ref{sec-walk}.
The proof of Theorem \ref{mainthm} is 
given at the end of Section \ref{sec-walk}.

\medskip
\noindent {\bf Acknowledgments.}
The author would like to thank Percy Deift 
for his interest and encouragement.

\section{Proof}\label{sec-walk}

Let $d_\lambda(N)$ be the number of semistandard Young tableaux 
of shape $\lambda$ with fillings taken from $\{1.2.\cdots,N\}$, 
and let $\ell(\lambda)$ be the number of rows of $\lambda$ 
(parts of $\lambda$, or the length of the first column).
From the bijection of path configurations and tableaux, the number of path 
$\pi\in\pathcon(N,k)$ satisfying $L_1(N,k)(\pi)\le l$ is equal to 
\begin{equation}\label{eq;walk3}   
  \sum_{\substack{\lambda\vdash k \\ \ell(\lambda)\le l}} d_\lambda(N).
\end{equation} 
In our analysis (see also \cite{BR1, BR2}), 
it turns out that in addition to 
the number of rows, the number of odd columns plays an important role 
in describing a tableau.
For a partition $\lambda$, we define 
$\lambda'$ to be the transpose of $\lambda$,
$f(\lambda)$ to be the number of odd row in $\lambda$,
and $|\lambda|$ to be the size of $\lambda$.
Let $b(N,j,m,l)$ be the number of semistandard Young tableaux
of size $2j+m$ with $m$ odd columns with at most $l$ columns
with fillings taken from $\{1,2,\cdots,N\}$ :
\begin{equation}\label{eq;walk4}
  b(N,j,m,l)= \sum_{\substack{\lambda\vdash 2j+m \\ f(\lambda')=m \\
\ell(\lambda)\le l}} d_\lambda(N).
\end{equation}
We use the notation $b(N,j,m,\infty)$ for the sum above without 
restriction on $\ell(\lambda)$.
Now we define a generating function 
\begin{equation}\label{eq-walk2}   
\begin{split}
  \phi(N,l,t,\beta) &:= 
(1-t^2)^{N(N-1)/2}(1-\beta t)^N 
\sum_{\ell(\lambda)\le l} t^{|\lambda|}\beta^{f(\lambda')}
d_\lambda(N)\\
&=(1-t^2)^{N(N-1)/2}(1-\beta t)^N
\sum_{j,m\ge 0} t^{2j}(\beta t)^m b(N,j,m,l),
\end{split}
\end{equation} 
where the sum in the first expression is taken 
over all the partitions $\lambda$ satisfying 
$\ell(\lambda)\le l$.

The starting point of our analysis is the following result of \cite{BR1}.
\begin{lem}\label{lem0}
Let $\phi(N,l,t,\beta)$ be defined as in \eqref{eq-walk2}. 
We have 
\begin{equation}\label{eq-walk3}   
  \phi(N,2l+1,t,\beta)
= (1-t^2)^{N(N-1)/2} \det(H_l), 
\end{equation} 
where $H_l=(h_{jk})_{0\le j,k<l}$ is the $l\times l$ Hankel determinant 
with 
\begin{equation}\label{eq.walk7}
  h_{jk}=h_{jk}(N,t)=\frac{2^{j+k+1}}{\pi} \int_{-1}^1 
x^{j+k} (1+t^2-2tx)^{-N}(1-x^2)^{1/2}dx.
\end{equation}
\end{lem}

\begin{rem}
Note that the right hand side of \eqref{eq-walk3} does not depend on $\beta$.
\end{rem}
\begin{proof}
This proof is in \cite{BR1} in a slightly different form.
Let $s_\lambda(x)$, $x=(x_1,x_2,\cdots)$, be the Schur function,  
and define $H(u;y)$ with $y=(y_1,y_2,\cdots)$ by 
\begin{equation}\label{eq;lem01}
  H(u;y)=\prod_j (1-uy_j)^{-1}.
\end{equation}
In (5.65) of \cite{BR1}, it is proved that 
\begin{equation}\label{eq;lem02}
  \sum_{\ell(\lambda)\le 2l+1} \beta^{f(\lambda')} s_\lambda(x)
= H(\beta;x) \Exp_{U\in Sp(2l)} \det(H(U;x)), 
\end{equation}
which is an identity as a formal power series in $x$.
But the combinatorial definition of the Schur function is 
(see, e.g. Chapter 7.10 of \cite{Stl}) 
\begin{equation}
  s_\lambda(x) =\sum_T x_1^{\alpha_1(T)} x_2^{\alpha_2(T)}\cdots,
\end{equation}
where the sum is over all semistandard Young tableaux $T$ of shape $\lambda$, 
and $\alpha_j(T)$ is the number of parts of $T$ equal to $j$ 
(type of $T$).
Since $\sum_j \alpha_j(T)=|\lambda|$, 
if we take the special case $x=(t,t,\cdots,t,0,0,\cdots)$ 
where the first $N$ elements are $t$ and the rest are $0$, 
then $s_\lambda(x)$ becomes 
\begin{equation}
  t^{|\lambda|}d_\lambda(N).
\end{equation}
Hence for this special choice of $x$ \eqref{eq;lem02} is now 
\begin{equation}\label{eq;lem05}
  \sum_{\ell(\lambda)\le 2l+1} \beta^{f(\lambda')} t^{|\lambda|}
d_\lambda(N)
= (1-\beta t)^{-N} \Exp_{U\in Sp(2l)} \det((1-tU)^{-N}).
\end{equation}

Now using Weyl's integration formula for $Sp(2l)$ (see, e.g. \cite{Simon}), 
the expectation in \eqref{eq;lem05} becomes 
\begin{equation}\label{eq;lem06}
  \Exp_{U\in Sp(2l)} \det((1-tU)^{-N})
= \frac{2^{l^2}}{l!(2\pi)^l}
\int_{[0,2\pi]^l} \prod_{1\le j<k\le l} (\cos\theta_j-\cos\theta_k)^2
\prod_{j=1}^{l} \sin^2\theta_j
(1+t^2-2t\cos\theta_j)^{-N}d\theta_j.
\end{equation}
Standard Vandermonde determinant manipulations yield 
\begin{equation}
  \frac{2^{l^2}}{(2\pi)^l} \det\biggl(
\int_0^{2\pi} \cos^{j+k}\theta \sin^2\theta(1+t^2-2t\cos\theta)^{-N}d\theta 
\biggr)_{0\le j,k<l},  
\end{equation}
which again after change of variables $x=\cos\theta$, 
is equal to 
\begin{equation}
  \frac{2^{l^2}}{\pi^l}\det\biggl(
\int_{-1}^1 U_j(x)U_k(x)(1+t^2-2tx)^{-N}(1-x^2)^{1/2}dx
\biggr)_{0\le j,k<l}
\end{equation}
where $U_j(x)=\frac{\sin((j+1)\theta)}{\sin\theta}$, $x=\cos\theta$, 
is the Chebyshev polynomial of the second kind.
Note that $U_j(x)=2^jx^j+\cdots$. 
Hence elementary row and column operations yield Lemma \ref{lem0}.
\end{proof}

Using this expression, we first obtain the asymptotic result 
for the generating function.
The limit is insensitive to $\beta$ since so is $\phi$.
\begin{prop}\label{prop1}
  Let $0<t<1$ and $\beta>0$ be fixed satisfying $0<\beta t<1$.
For each $l$ and $N$, define $x\in\R$ by 
\begin{equation}\label{eq.walk16}   
  x=\frac{l-\eta(t)N}{\rho(t)N^{1/3}},
\end{equation} 
where $\eta(t)$ and $\rho(t)$ are defined in \eqref{eq;walk1}
Then there exits a positive constant $M_0$ such that for $M>M_0$, 
there are constants $C, c>0$, independent of $M$, and $C(M)$ 
which may depend on $M$ so that 
\begin{equation}\label{eq;walk7}   
  |\phi(N,l,t,\beta)- F_1(x)| \le 
\frac{C(M)}{l^{1/3}}+ Ce^{-cM^{3/2}}, 
\qquad -M<x<M.
\end{equation} 
Also we have 
\begin{eqnarray}
\label{walk16}
  0\le1-\phi(N,l,t,\beta)\le Ce^{-cx^{3/2}}, \qquad &x\ge M_0,\\
\label{walk17}
  0\le \phi(N,l,t,\beta)\le Ce^{-c|x|^3}, \qquad &x\le -M_0.
\end{eqnarray}
\end{prop}

\begin{proof}
  It is enough to consider the limit for $\phi(N,2l+1,t,\beta)$ 
since from the definition \eqref{eq-walk2}, $\phi$ is monotone
increasing in $l$.
First we related the determinant in \eqref{eq-walk3} with certain quantities 
of orthogonal polynomials on the circle.

Let $p_j(x)=x^j+\cdots$ be the $j^{\text{th}}$ monic orthogonal 
polynomial with respect to the weight $w(x)dx=(1+t^2-2tx)^{-N}(1-x^2)^{1/2}dx$ 
on the interval $(-1,1)$, and let $C_j$ be the norm of $p_j$ : 
\begin{equation}
  \int_{-1}^1 p_j(x)p_k(x) w(x)dx= C_j\delta_{jk}.
\end{equation}
It is a well known result of orthogonal polynomial theory 
(see, e.g. \cite{Szego}) 
that $C_j=\det(\tilde{H}_{j+1})/\det(\tilde{H}_j)$, 
where $\tilde{H}_l=(\tilde{h}_{jk})_{0\le j,k<l}$ 
with 
\begin{equation}
  \tilde{h}_{jk}=\int_{-1}^1 x^{j+k}(1+t^2-2tx)^{-N}
(1-x^2)^{1/2}dx.
\end{equation}
which is equal to $\frac{\pi}{2^{j+k+1}}h_{jk}$ (recall \eqref{eq.walk7}).
Hence $\det(\tilde{H}_j)=\frac{\pi^j}{2^{j^2}}\det(H_j)$.
Since the Szeg\"o strong limit theorem for Hankel determinants 
(see, e.g. \cite{Jo1}) implies that 
$\lim_{l\to\infty}\det(H_l)= (1-t^2)^{-N(N-1)/2}$ for fixed $N$,  
we have 
\begin{equation}\label{eq.walk23}
  \phi(N,2l+1,t,\beta)=\lim_{k\to\infty} 
\prod_{j=l}^{k}\frac{\det(H_j)}{\det(H_{j+1})}
=\prod_{j=l}^{\infty} \frac{\pi}{2^{2j+1}}C_j^{-1}.
\end{equation}

Now we use the relation between orthogonal polynomials on the unit circle 
and those on the interval $(-1,1)$.
Let $\pi_j(z)=z^j+\cdots$ be the $j^{\text{th}}$ monic orthogonal polynomials 
on the unit circle $\{|z|=1\}$ with respect to the weight 
\begin{equation}
  \varphi(z)\frac{dz}{2\pi i}=(1-tz)^{-N}(1-tz^{-1})^{-N}\frac{dz}{2\pi iz},
\end{equation}
and let 
$N_j$ be the norm of $\pi_j(z)$ : 
\begin{equation}
  \int_{|z|=1} \pi_j(z)
\overline{\pi_k(z)}\varphi(z) \frac{dz}{2\pi iz}
= N_j \delta_{jk}.
\end{equation}
There is a simple relation between orthogonal polynomials $p_j$ on 
the unit circle and orthogonal polynomials $\pi_j$ on the interval 
(see the forth equation of (11.5.2) in \cite{Szego}) :
\begin{equation}
  C_j^{-1/2}p_j(x) = \sqrt{\frac{2}{\pi}} \bigl(1+\pi_{2j+2}(0)\bigr)^{-1/2}
N_{2j+1}^{-1/2} 
\frac{z^{-j}\pi_{2j+1}(z)-z^{j}\pi_{2j+1}(z^{-1})}{z-z^{-1}}, 
\quad x=\frac12(z+z^{-1}).
\end{equation} 
Especially comparing the coefficient of the leading term $x^j$, 
we have the relation 
\begin{equation}\label{eq.walk26}
   C_j=\frac{\pi}{2^{2j+1}} \bigl(1+\pi_{2j+2}(0)\bigr)N_{2j+1}.
\end{equation}
But we also have 
$(1-\pi_{n+1}(0)^2)N_n=N_{n+1}$ (see (11.3.6) in \cite{Szego}). 
Hence \eqref{eq.walk26} is equal to 
\begin{equation}
   C_j=\frac{\pi}{2^{2j+1}} \bigl(1-\pi_{2j+2}(0)\bigr)^{-1}N_{2j+2}.
\end{equation}
Therefore \eqref{eq.walk23} becomes 
\begin{equation}
  \phi(N,2l+1,t,\beta) = 
\prod_{j=l}^\infty \bigl(1-\pi_{2j+2}(0)\bigr)N_{2j+2}^{-1}.
\end{equation}
This argument appeared in Corollary 2.7 of \cite{BR1}.
Now using Proposition \ref{prop-rhp}, computations similar to 
Lemma 7.1 of \cite{BDJ} (also Corollary 7.2 and 
Corollary 7.6 of \cite{BR2}) 
yield the result.
\end{proof}

Interpreting the notation $b(N,j,m,\infty)$ as the sum without 
constraints on $\ell(\lambda)$ in \eqref{eq;walk4}, 
the number of path configuration in $\pathcon(N,k)$ is equal to 
\begin{equation}\label{eq;walk8}
   |\pathcon(N,k)| = \sum_{2j+m=k} b(N,j,m,\infty), 
\end{equation}
and the probability of interest that $L_1(N,k)\le l$ for 
$\pi\in\pathcon(N,k)$ is equal to 
\begin{equation}\label{eq;walk9}
   \frac1{|\pathcon(N,k)|}
\sum_{2j+m=k} b(N,j,m,l)
= \frac1{|\pathcon(N,k)|}
\sum_{2j+m=k} p(N,j,m,l)b(N,j,m,\infty),
\end{equation}
where 
\begin{equation}\label{eq.walk30}
 p(N,j,m,l)=\frac{b(N,j,m,l)}{b(N,j,m,\infty)}.
\end{equation}

For fixed $N$, as $l\to\infty$, 
the Szeg\"o strong limit theorem for Hankel determinants 
(see, e.g. \cite{Jo1}) implies that 
\eqref{eq-walk3} becomes $1$.
Thus we have the identity 
\begin{equation}\label{eq;walk10}   
  \sum_{j,m\ge 0} t^{2j}(\beta t)^{m} b(N,j,m,\infty)
= (1-t^2)^{-N(N-1)/2}(1-\beta t)^{-N}.
\end{equation} 
By taking Taylor expansion of the right hand side in $t$ and $\beta$, 
we obtain 
\begin{equation}\label{eq-walk8}   
  b(N,j,m,\infty)
= \binom{\frac{N(N-1)}2+j-1}{j}\binom{N+m-1}{m}.
\end{equation} 
There is a more direct way to see this.
See the remark after Lemma \ref{lem3} below.
Now from \eqref{eq;walk8}, the total number of paths 
in $\pathcon(N,k)$ is equal to 
\begin{equation}\label{eq;walk12}   
   |\pathcon(N,k)|=
\sum_{2j+m=k} \binom{\frac{N(N-1)}2+j-1}{j}\binom{N+m-1}{m}.
\end{equation} 
Now it is straightforward to obtain the following result on the number 
of all paths.

\begin{lem}\label{lem1}
 Let $0<t<1$
and let 
\begin{equation}\label{eq;walk13}
  k= \bigl[ \frac{t^2}{1-t^2}N^2+o(N^{4/3})\bigr]. 
\end{equation}
As $N\to\infty$, we have 
\begin{equation}\label{eq.walk33}
  |\pathcon(N,k)| = 
\frac{\exp\bigl\{-\frac{N^2t^2}{1-t^2}\log t 
-\mu N\log t- \frac14\bigl(\frac{\mu(1-t^2)}{t}-1\bigr)^2\bigr\}}
{\sqrt\pi tN(1-t)^{N}(1-t^2)^{N(N-1)/2-1}}
\bigl(1+o(1)\bigr)
\end{equation}
where the term $o(1)$ vanishes as $N\to\infty$, 
and $\mu$ is defined by 
\begin{equation}\label{eq.walk34}
  \mu:= \frac{k}{N}-\frac{t^2}{1-t^2}N
\end{equation}
which is of order $o(N^{1/3})$ from \eqref{eq;walk13}.
Moreover, the main contribution to the sum comes from 
$|m - \frac{t}{1-t}N| \le N^{1/2+\epsilon/2}$ for some $0<\epsilon<\frac13$ ; 
precisely, there is a constant $c>0$ such that 
for any $0<\epsilon<\frac13$, we have 
\begin{equation}\label{eq.walk37}
  |\pathcon(N,k)|=
\biggl[ \sum_{\substack{|m-\frac{t}{1-t}N|\le N^{1/2+\epsilon/2} \\ 
k-m \text{ is even}}} b(N,\frac{k-m}{2},m,\infty) \biggr]
\bigl(1+O(e^{-cN^{\epsilon}})\bigr).
\end{equation}
\end{lem}

\begin{proof}
From \eqref{eq;walk12}, we have 
\begin{equation}
  |\pathcon(N,k)|=\sum_{m=0}^{[\frac{k}2]} a(m), 
\qquad a(m)= \binom{\frac{N(N-1)}{2}+\frac{k-m}2-1}{\frac{k-m}2}
\binom{N+m-1}{m}.
\end{equation}
The ratio of $a(m)$ is 
\begin{equation}\label{eq.walk38}
  \frac{a(m+2)}{a(m)}=\frac{(N+m+1)(N+m)(k-m)}{(m+2)(m+1)(N(N-1)+k-m-2)}.
\end{equation}
One can directly check that under the condition \eqref{eq;walk13}, 
the above ratio is decreasing in $m$, and becomes closest to $1$ at 
\begin{equation}
  m_c= \bigl[\frac{t}{1-t}N + o(N^{1/3})\bigr].
\end{equation}
Hence $a(m)$ is unimodal: it is increasing for $m<m_c$ 
and is decreasing for $m>m_c$.
Now consider the neighborhood $\mathcal N$ 
of $m_c$ of size $N^{1/2+\epsilon/2}$ 
for some fixed $0<\epsilon<\frac13$.
For $m$ in $\mathcal N$, set 
\begin{equation}\label{eq.walk39}
  m= \frac{t}{1-t}N+x, \qquad |x|\le N^{1/2+\epsilon/2}.
\end{equation}
For any $M,x>0$, Stirling's formula yields 
\begin{equation}\label{2.41}
  (M+x)!=\sqrt{2\pi M}M^{M+x}e^{-M+\frac{x^2}{2M}}\
\biggl(1+O\bigl(\frac{x}{M}\bigr)+O\bigl(\frac{x^3}{M^2}\bigl)\biggr).
\end{equation}
Using \eqref{eq;walk13}, \eqref{eq.walk39} and \eqref{2.41}, 
we have for $m$ in $\mathcal N$, 
\begin{equation}   
  \binom{N+m-1}{m}=\frac{(\frac{1}{1-t}N+x-1)!}{(N-1)!(\frac{t}{1-t}N+x)!}
=\frac{\exp\bigl\{-(\frac{tN}{1-t}+x)\log t-\frac{(1-t)^2x^2}{2tN}\bigr\}}
{\sqrt{2\pi tN}(1-t)^{N-1}}
\biggl(1+O(N^{-\frac12+\frac32\epsilon})\biggr)
\end{equation} 
and 
\begin{equation}   
\begin{split}
  &\binom{\frac{N(N-1)}{2}+\frac{k-m}2-1}{\frac{k-m}2}
=\frac{\bigl(\frac{N(N-1)}{2(1-t^2)}+(\mu-\frac{t}{1-t^2})\frac{N}2
-\frac{x}2-1\bigr)!}
{\bigl(\frac{N(N-1)}{2}-1\bigr)!
\bigl(\frac{t^2N(N-1)}{2(1-t^2)}+(\mu-\frac{t}{1-t^2})\frac{N}2
-\frac{x}2\bigr)!}\\
&\qquad\qquad 
=\frac{\exp\bigl\{ -(\frac{t^2}{1-t^2}N^2+(\mu-\frac{t}{1-t})N-x)\log t 
-\frac{1}{4}(\frac{(1-t^2)\mu}{t}-1)^2  \bigr\}}
{\sqrt{\pi}tN(1-t^2)^{N(N-1)/2-1}} \bigl(1+o(1)\bigr).
\end{split}
\end{equation} 
Thus we have 
\begin{equation}\label{eq.walk40}
  a(m) = 
\frac{\exp\bigl\{-\frac{N^2t^2}{1-t^2}\log t 
-N\mu\log t- \frac14\bigl(\frac{\mu(1-t^2)}{t}-1\bigr)^2
-\frac{(1-t)^2x^2}{2tN}\bigr\}}
{\sqrt{2}\pi t^{3/2}N^{3/2}(1-t)^{N}(1-t^2)^{N(N-1)/2-1}}
\bigl(1+o(1)\bigr).
\end{equation}
Let 
\begin{equation}
  |\pathcon(N,k)| =\sum_{*} a(m)+\sum_{**} a(m), 
\end{equation}
where $*$ denotes the set $\mathcal N$ of $m$ satisfying \eqref{eq.walk39} 
and $**$ denotes the rest of the range of $m$.
From \eqref{eq.walk40}, the first sum over $*$ is equal 
to the right hand side of \eqref{eq.walk33}.
Also from the unimodality, $a(m)$ in $**$ is less than or equal to 
the largest of $a(m_+)$ and $a(m_-)$ where 
$m_\pm=[\frac{t}{1-t}N \pm N^{1/2+\epsilon}]$.
The number of summand in $**$ is of order $N^2$.
Hence using \eqref{eq.walk40} again for $m_\pm$, 
if we take $c=\frac{(1-t)^2}{4t}$, for large $N$, 
we obtain 
\begin{equation}
  \sum_{**} a(m)= \bigl(\sum_{*} a(m)\bigr) e^{-cN^{\epsilon}}, 
\end{equation}
which establishes the proof.
\end{proof}

Now we rewrite \eqref{eq-walk2} as 
\begin{equation}
  \phi(N,l,t,\beta)
= (1-t^2)^{N(N-1)/2}(1-\beta t)^N
\sum_{j,m\ge 0} t^{2j}(\beta t)^m b(N,j,m,\infty) p(N,j,m,l).
\end{equation}
The asymptotics of $\phi(N,l,t,\beta)$ and  $p(N,j,m,l)$ 
are related as follows.
\begin{lem}\label{lem-dep}
  For any $d>0$, there are constants $C_0, c_0>0$ 
such that 
for all $l\ge 0$, 
\begin{equation}\label{eq.walk46}
  p(N,\mu_+(N),\nu_+(N),l)-\frac{C_0}{N^d}\le 
\phi(N,l,t,\beta)\le p(N,\mu_-(N),\nu_-(N),l)+\frac{C_0}{N^d}
\end{equation}
where
\begin{eqnarray}
\label{eq.walk47}
  \mu_\pm(N)&=&\biggl[\frac{t^2}{2(1-t^2)}N^2 \pm c_0N\sqrt{\log N}\biggr],\\
\label{eq.walk48}
  \nu_\pm(N)&=&\biggl[\frac{\beta t}{1-\beta t} N \pm c_0\sqrt{N\log N}\biggr].
\end{eqnarray}
\end{lem}
The proof follows by using the following 
Lemma twice for $j$ and $m$ indices together with 
the Lemma \ref{lem3}. (Recall the \eqref{eq-walk8}).

\begin{lem}\label{lem2}
For a sequence $\{q_j\}_{j\ge 0}$, we define the following 
generating function
\begin{equation}\label{eq46}
  G(N)=(1-a)^N\sum_{j=0}^\infty a^j\binom{N+j-1}{j} q_j, 
\qquad 0<a<1,\ \ N=1,2,\cdots.
\end{equation}
  For each $d>0$, there are constants $C_1,c_1\ge 0$ such that 
for any sequence $\{q_j\}_{j\ge 0}$ satisfying 
(i) $q_j\ge q_{j+1}$ and (ii) $0\le q_j\le 1$, 
\begin{equation}
   q_{N^{**}}-\frac{C_1}{N^d}\le G(N)\le q_{N^*}+\frac{C_1}{N^d}, 
\qquad N\ge 1,
\end{equation}
where 
\begin{eqnarray}
  N^*&=&\frac{a}{1-a}N-c_1\sqrt{N\log{N}},\\
  N^{**}&=&\frac{a}{1-a}N+c_1\sqrt{N\log{N}}.
\end{eqnarray}
\end{lem}

\begin{proof}
  This proof is parallel to that of the de-Poissonization 
lemma in \cite{Jo2}. 
We have 
\begin{equation}
  G(N)=\sum_{j=0}^\infty f_jq_j, \qquad f_j=(1-a)^Na^j\binom{N+j-1}{j}.
\end{equation}
Stirling's formula yields for $n,m\ge 1$, 
\begin{equation}
  \binom{m+n-1}{m}\le C\exp\biggl\{ m\bigl[(1+\frac{n}{m}\log(1+\frac{n}{m})
-\frac{n}{m}\log\frac{m}{n}\bigr] \biggr\},
\end{equation}
with some constant $C$.
Thus we have 
\begin{equation}
  f_j\le C\exp\{Nh(j/N)\}, 
\qquad h(x)=(1+x)\log(1+x)-x\log x +x\log a+\log(1-a).
\end{equation}
One can directly check the following estimates of $h$ :
\begin{eqnarray}
\label{eq53}
  h(x) &\le& -\frac{(1-a)^2}{2a}\bigl(x-\frac{a}{1-a}\bigr)^2, 
\qquad \qquad 0\le x\le \frac{2a}{1-a},\\
\label{eq54}
  h(x) &\le& -\bigl[ \log 2-\frac{1+a}{2a}\log(1+a)\bigr] x, 
\qquad x\ge \frac{2a}{1-a}.
\end{eqnarray}
We take a constant $c_1>0$ satisfying $\frac{(1-a)^2}{2a}c_1^2-1\ge d$.
From \eqref{eq53} and the condition (ii), we have 
\begin{equation}
  \sum_{j\le N^*} f_jq_j \le \frac{a}{1-a}CN 
e^{-\frac{(1-a)^2}{2a}c^2\log N}
\le \frac{C}{N^d}, 
\end{equation}
with a new constant $C$.
Similarly, 
\begin{equation}
  \sum_{N^{**}\le j\le \frac{2a}{1-a}N} f_jq_j \le \frac{C}{N^d}.
\end{equation}
Also using \eqref{eq54}, we have 
\begin{equation}
  \sum_{j\ge \frac{2a}{1-a}N} f_jq_j\le C'e^{-c'N}
\end{equation}
for some constants $c',C'$.
Thus we have 
\begin{equation}\label{eq58}
 \biggl| G(N)-\sum_{N^*\le j\le N^{**}}f_jq_j \biggr| \le \frac{C}{N^d},
\end{equation}
with a possibly different constant $C$.
Now from the monotonicity condition (i), we have  
\begin{equation}
  \sum_{N^*\le j\le N^{**}} f_jq_j
\le \biggl( \sum_{N^*\le j\le N^{**}} f_j \biggr) q_{N^*}
\le q_{N^*},
\end{equation}
and 
\begin{equation}
  \sum_{N^*\le j\le N^{**}} f_jq_j
\ge \biggl(  \sum_{N^*\le j\le N^{**}} f_j\biggr) q_{N^{**}}
\ge \bigl(1-\frac{C}{N^d}\bigr) q_{N^{**}}
\ge q_{N^{**}} -\frac{C}{N^d}, 
\end{equation}
using the equality \eqref{eq58} for the second equality.
Thus we obtained the desired result.
\end{proof}

To use the above Lemma to $\phi$, we need monotonicity in $l$.
It is more convenient now to view semistandard Young tableaux (SSYT) 
as generalized permutations.
A two-rowed array
\begin{equation}
  \pi=\begin{pmatrix}
i_1&\cdots&i_k\\ j_1&\cdots&j_k
\end{pmatrix}
\end{equation}
is called a generalized permutation if either $i_r<i_{r+1}$
or $i_r=i_{r+1}$, $j_r\le j_{r+1}$.
Suppose the elements in the upper row of $\pi$ come from 
$\{1,2,\cdots,M\}$ and the elements in the bottom row 
come from $\{1,2,\cdots,N\}$.
One can represent a generalized permutation as a $M\times N$ matrix 
$(a_{ik})$ where $a_{ik}$ is the number of times when 
$\binom{i}{k}$ occurs in $\pi$.
For example, the generalized permutation
\begin{equation}\label{eq;walk68}
\begin{pmatrix}
 1&1&1&2&2&2&2&3&3\\ 1&3&3&2&2&2&4&3&4
\end{pmatrix}
\end{equation}
corresponds to
\begin{equation}
\begin{pmatrix}
  1&0&2&0\\ 0&3&0&1\\ 0&0&1&1
\end{pmatrix}.
\end{equation}
In the proof of the Lemma below, we regard a generalized permutation 
as a $M\times N$ square board with stacks of $a_{ik}$ 
balls in each position $(i,k)$. 
We denote by $L(\pi)$ the length of the longest strictly decreasing
subsequence of $\pi$.
In the example \eqref{eq;walk68}, $L(\pi)=2$.

Let $M_{j,m}$ be the set of $N\times N$ matrices $\pi=(a_{ik})$
which is symmetric $a_{ik}=a_{ki}$, and satisfies 
$\sum_{i=1}^N a_{ii}= m$ and $\sum_{1\le i<k\le N} a_{ik}= j$.
This is a certain subset of the set of generalized permutations.
The celebrated Robinson-Schensted-Knuth correspondence 
\cite{Knuth:correspondence}
establishes a bijection between $M_{j,m}$ 
and the set of SSYT of size $2j+m$ 
with $m$ odd columns with fillings taken from 
$\{1,2,\cdots, N\}$.
Moreover, under this bijection, 
$L(\pi)$ for $\pi\in M_{j,m}$ (viewed as a generalized permutation) 
is equal to the number of rows of the corresponding SSYT.

With this preliminary, we can prove the following.

\begin{lem}[monotonicity]\label{lem3}
  For any $j,m\ge 0$, we have 
\begin{equation}
  p(N,j+1,m,l)\le p(N,j,m,l), \qquad p(N,j,m+1,l)\le p(N,j,m,l).
\end{equation}
\end{lem}

\begin{proof}
 We first consider the second inequality.
From \eqref{eq.walk30}, we need to show that 
\begin{equation}
  (m+1)b(N,j,m+1,l)\le (N+m)b(N,j,m,l).
\end{equation}
By the definition \eqref{eq;walk4} and 
the Robinson-Schensted-Knuth correspondence, 
$b(N,j,m,l)$ is equal to the number of $\pi\in M_{j,m}$ 
satisfying $L(\pi)\le l$.

Consider all possible distinct (strict) upper triangular parts 
of elements in $M_{j,m}$.
It is equal to putting $j$ identical balls into $N(N-1)/2$ boxes ; 
$K=\binom{N(N-1)/2+j-1}{j}$ distinct ways.
Hence we have a disjoint union $M_{j,m}=\cup_{i=1}^{K} S_{i,m}$ 
where each $S_{i,m}$ consists of $\pi\in M_{j,m}$ with same upper 
triangular part, and elements in $S_{i,m}$ and $S_{i',m}$ 
have different upper triangular parts when $i\neq i'$.
Similarly, $M_{j,m+1}=\cup_{i=1}^{K} S'_{i,m+1}$
where $\sigma\in S'_{i,m+1}$ has the upper triangular part same as 
that of $\pi\in S_{i,m}$.

Now for each $\pi=(a_{rs})\in S_{i,m}$, we generate 
$N+m$ elements in $S'_{i,m+1}$ as follows. 
For $1\le r\le N$, assign $a_{rr}+1$ identical $\pi'=(a'_{kl})$ 
such that $a'_{rr}=a_{rr}+1$ and $a'_{kl}=a_{kl}$ 
for $(k,l)\neq (r,r)$.
(One can think this as adding a new ball in an array of $a_{rr}$ 
balls ; there are $a_{rr}+1$ ways.)
Since $\sum_{1\le r\le N} a_{rr}+1= m+N$, 
there result $(m+N)|S_{i,m}|$ (many identical) elements of $S_{i,m+1}$.
Note that under this assignment, 
\begin{equation}\label{eq;walk70}
  L(\pi')\ge L(\pi).
\end{equation}
Now fix $\sigma=(b_{kl})\in S_{i,m+1}$.
Since there are $m+1$ diagonal entries, there are exactly 
$m+1$ (many identical) elements of $S_{i,m}$ from which 
$\sigma$ is generated under the above assignment.
(Considering each entry as a ball, each $m+1$ balls on the diagonal 
can be a newly added one.)
Thus we have the identity $(m+N)|S_{i,m}|=(m+1)|S_{i,m+1}|$.
Furthermore, from the remark regarding \eqref{eq;walk70}, 
we have $(m+1)|R(i,m+1,l)|\le (m+N)|R_{i,m,l}|$, 
where $R_{i,m,l}$ is the subset of $\pi\in S_{i,m}$ satisfying 
$L(\pi)\le l$.
Therefore the second inequality in the Lemma is obtained.

The first inequality follows from a similar argument.
\end{proof}

\begin{rem}
As mentioned before, using the generalized permutation interpretation 
of SSYT, we can see \eqref{eq-walk8} directly.
From non-negative integer matrix representation of generalized permutations, 
$b(N,j,m,\infty)$ is the number of $N\times N$ matrices $(a_{rs})$ 
with non-negative integer entries such that 
$\sum_{r=1}^N a_{rr}=m$ and 
$\sum_{1\le r<s\le N}a_{rs}=j$.
It is equivalent to placing $m$ identical balls of color 1 into $N$ boxes 
and $j$ identical balls of color 2 into $N(N-1)/2$ boxes.
Therefore we obtain \eqref{eq-walk8}.
\end{rem}

Now we give the proof of the theorem.
\begin{proof}[proof of Theorem \ref{mainthm}]
In \eqref{eq;walk9}, we have 
\begin{equation}
  \Prob(L_1(N,k)\le l) = \frac{1}{|P(N,k)|}
\sum_{2j+m=k} p(N,j,m,l)b(N,j,m,\infty).
\end{equation}
where $P(N,k)$, $p(N,j,m,l)$ and $b(N,j,m,\infty)$ 
are given in \eqref{eq;walk12}, \eqref{eq-walk8} and \eqref{eq.walk30}, 
and $b(N,j,m,l)$ is given in \eqref{eq;walk4}.
We split the above sum into two pieces. 
One part is the sum over (1) $|m-\frac{t}{1-t}N|\le N^{1/2+\epsilon/2}$
and (2) the rest, where $0<\epsilon<\frac13$ is fixed.
Then since $0\le p(N,j,m,l)\le 1$, we have from \eqref{eq.walk37}
\begin{equation}
  \frac{1}{|P(N,k)|} \sum_{(2)} p(N,j,m,l) b(N,j,m,\infty)
= O(e^{-cN^{\epsilon}})
\end{equation}
for some $c>0$.

We use Lemma \ref{lem-dep} to estimate $p(N,j,m,l)$ for $(j,m)$ in (1).
Set
\begin{eqnarray}
  \tilde{t}^2=\frac{k-m-2c_0N\sqrt{\log N}}{N^2+ k-m-2c_0N\sqrt{\log N}},
\qquad
\beta=\frac{m-c_0\sqrt{N\log N}}{N+ m-c_0\sqrt{N\log N}}
\tilde{t}^{-1} ,
\end{eqnarray}
where $k$ satisfies the condition in \eqref{eq;walk2},
and we take $\tilde{t}>0$.
For $(j,m)$ in (1), they satisfy
\begin{eqnarray}\label{eq.walk74}
\tilde{t}=t+o(N^{-2/3}),        
\qquad \beta=1+O(N^{-1/2+\epsilon/2}).
\end{eqnarray}
The first inequality of \eqref{eq.walk46} yields
$p(N,j,m,l)\le \phi(N,l,\tilde{t},\beta)+C_0N^{-d}$.
Set $\tilde{x}$ by \eqref{eq.walk16} where $t$ is replaced by $\tilde{t}$
and $l$ is given by $l=[\eta(t)N+x\rho(t)N^{1/3}]$.
Let $M?M_0$ satisfies $-M<2x<M$ where $M_0$ is given in 
Proposition \ref{prop1}.
From \eqref{eq.walk74}, we have 
\begin{equation}\label{2.76}
  \tilde{x}=x+o(1),
\end{equation}
Proposition \ref{prop1} implies that 
\begin{equation}
|\phi(N,l,\tilde{t},\beta) - F_1(\tilde{x})|\le 
\frac{C(M)}{l^{1/3}}+Ce^{-cM^{3/2}}, 
\qquad l=[\eta(t)N+x\rho(t)N^{1/3}]
\end{equation}
for large $N$.
Since $F_1'(x)=-\frac12(u(x)+v(x))F_1(x)$ is bounded for $x\in\R$, 
we have from \eqref{2.76} that 
$F_1(\tilde{x})=F_1(x)+ o(1)$.     
Thus we have for large $N$, 
\begin{equation}
  p(N,j,m,l)\le F_1(x) + o(1),
\end{equation}
where $o(1)$ term is independent of $(j,m)$ in (1)
and vanishes as $N\to\infty$.
Thus using Lemma \ref{lem1}, we have for large $N$, 
\begin{equation}
\begin{split}
  \Prob(L_1(N,k)\le l) &\le \frac{1}{|P(N,k)|} 
\sum_{(1)} (F_1(x)+o(1)) b(N,j,m,\infty) + O(e^{-cN^{2\epsilon}})\\
&= F_1(x) + o(1), 
\qquad l=[\eta(t)N+x\rho(t)N^{1/3}]
\end{split}
\end{equation}

Similarly, we obtain the lower bound using the second inequality 
of \eqref{eq.walk46}.
Thus we proved \eqref{eq;walk2}.

The convergence of moments is also similar 
using \eqref{walk16} and \eqref{walk17}
(cf. Section 8 of \cite{BR2}).
\end{proof}

\section{Asymptotics of orthogonal polynomials}\label{sec-rhp}

This section is devoted to asymptotics of orthogonal polynomials 
used in the proof of Proposition \ref{prop1}. 
The key ingredient is the equilibrium measure 
(see \cite{DKMVZ2, DKMVZ3, BDJ}).

On the unit circle, the equilibrium measure 
$d\mu_V(z)=\psi(\theta)\frac{d\theta}{2\pi}$ 
for $V(z)$ and its support are uniquely determined by the following 
Euler-Lagrange variational conditions :
\begin{equation}\label{eq-varcon}
  \begin{split}
  &\text{there exits a real constant $l$ such that},\\
  &\quad 2\int_{\Sigma} \log{|z-s|}d\mu_V(s) - V(z)+l=0 \ \ \text{for $z\in \bar
{J}$},\\
  &\quad 2\int_{\Sigma} \log{|z-s|}d\mu_V(s) - V(z)+l \le 0 \ \
\text{for $z\in\Sigma -\bar{J}$}.
  \end{split}
\end{equation}

\begin{lem}\label{as-lem1}
Let $\gamma\ge 0$ and $0<t<1$, and let 
\begin{eqnarray}
V(z)&=&\gamma\log(1-t z)(1-t z^{-1}).
\end{eqnarray}
Then their equilibrium measure $\psi(\theta)d\theta/2\pi$ 
is given as follows.
\begin{itemize}
\item When $0\le\gamma\le \frac{1+t}{2t}$, we have 
$J=\Sigma$, $l=0$, and 
\begin{eqnarray}
  \psi(\theta) 
&=& 1-\gamma+\frac{\gamma (1-t^2)}{1+t^2-2t\cos\theta}\\
&=& 1+\gamma t\biggl( \frac{z}{1-tz}+\frac{z^{-1}}{1-tz^{-1}} \biggr), 
\qquad z=e^{i\theta}.
\end{eqnarray}
\item When $\gamma>\frac{1+t}{2t}$, 
$J=\{e^{i\theta} : |\theta|\le\theta_c\}$, where 
$\sin^2\frac{\theta_c}2=\frac{(1-t)^2(2\gamma-1)}{4t(\gamma-1)^2}$, 
$0<\theta_c<\pi$, 
or 
\begin{equation}
  |\xi-t|=\frac{\gamma(1-t)}{\gamma-1}, \quad \xi=e^{i\theta_c}, 
\end{equation}
and 
\begin{equation}
  l= 2\gamma\log\frac{(2\gamma-1)(1-t)}{2(\gamma-1)} 
-\log\frac{(2\gamma-1)(1-t)^2}{4t(\gamma-1)^2},
\end{equation}
and finally 
\begin{eqnarray}
  \psi(\theta)
&=& \frac{4(\gamma-1)\cos
\frac{\theta}{2}}{\frac{(1-t)^2}{t}+4\sin^2\frac{\theta}2}
\sqrt{\sin^2\frac{\theta_c}2-\sin^2\frac{\theta}2} \\
&=& \frac{(\gamma-1)(1+z)}{(z-t)(z-t^{-1})}\sqrt{(z-\xi)(z-\xi^{-1})_+}, 
\qquad z=e^{i\theta},
\end{eqnarray}
where $\sqrt{(z-\xi)(z-\xi^{-1})_+}$ denotes the limit 
of $z$ from inside the unit circle.
And in this case, the inequality of the second variational condition 
in \eqref{eq-varcon} is strict for $z\in\Sigma\setminus\overline{J}$.
\end{itemize}
\end{lem}

\begin{proof}
The proof given here is similar to the proof of Lemma 4.3 in \cite{BDJ}, whose 
main ingredient is the following results of Lemma 4.2 in \cite{BDJ}.
Let $d\mu(s)=u(\theta)d\theta$ be an absolutely continuous probability 
measure on the unit circle $\Sigma$ and $u(\theta)=u(-\theta)$.
Define 
\begin{equation}\label{eq-gdef}
  g(z)=\int_{\Sigma} \log(z-s) d\mu(s), 
\end{equation}
where for fixed $s=e^{i\theta_0}\in\Sigma$, $\log(z-s)$ is defined to be 
analytic in $\C\setminus\bigl((-\infty,-1]\cup
\{e^{i\theta} : -\pi\le\theta\le\theta_0\}\bigr)$, 
and $\log(z-s)\sim \log z$ for $z\to +\infty$ with $z\in\R$.
Then for $z=e^{i\phi}\in\Sigma$, we have 
\begin{eqnarray}
\label{eq-g1}
g_+(z)+g_-(z)&=&2\int_\Sigma \log|z-s|d\mu(s) + i(\phi+\pi), \\
\label{eq-g2}
g_+(z)-g_-(z)&=&2\pi i\int_{\phi}^{\pi} u(\theta)d\theta.
\end{eqnarray}
Also evenness of $u$ yields $g(0)=\pi i$.
\begin{itemize}
\item When $0\le\gamma\le \frac{1+t}{2t}$ : 
By residue calculation, it is easy to check 
$\int_{-\pi}^\pi\psi(\theta)\frac{d\theta}{2\pi}=1$.
Define $g(z)=\int_{-\pi}^\pi \log(z-e^{i\theta}) \psi(\theta) 
\frac{d\theta}{2\pi}$ as in \eqref{eq-gdef}.
Then by direct residue calculation, we obtain 
\begin{equation}   
g'(z)=
 \begin{cases}
   \frac{-\gamma t}{1-tz}, \qquad |z|<1,\\
   \frac{\gamma tz^{-2}}{1-tz^{-1}}+\frac{1}{z}, \qquad |z|>1.
 \end{cases}
\end{equation} 
Since $g(0)=\pi i$ and $g(z)\sim\log(z)$ as $z\to\infty$, 
we have 
\begin{equation}   
g(z)= \begin{cases}
 \gamma\log(1-tz)+\pi i, \qquad |z|<1,\\
 \gamma\log(1-tz^{-1})+\log z, \qquad |z|>1.
\end{cases}
\end{equation} 
Thus, from \eqref{eq-g1}, the variational condition \eqref{eq-varcon} 
is satisfied with $J=\Sigma$ and $l=0$.
\item 
When $\gamma>\frac{1+t}{2t}$ :
Set $\beta(z)=\sqrt{(z-\xi)(z-\xi^{-1})}$ which is analytic in 
$\C\setminus \overline{J}$ and $\beta(z)\sim z$ as $z\to +\infty$ 
with $z\in\R$.
Then we have  
\begin{equation}\label{eq-rhp13}
  \beta(0)=-1, \quad \beta(t)=-|\xi-t|, \quad 
\beta(t^{-1})=\frac1{t}|\xi-t|.
\end{equation}
First, direct residue calculation using \eqref{eq-rhp13} shows that 
$\int_{-\theta_c}^{\theta_c}\psi(\theta)\frac{d\theta}{2\pi}=1$, 
hence $\psi(\theta)d\theta/(2\pi)$ is a probability measure.
Now we define $g(z)=\int_{-\theta_c}^{\theta_c}\log(z-e^{i\theta}) 
\psi(\theta)\frac{d\theta}{2\pi}$ as before.
Using \eqref{eq-rhp13} again, residue calculations yield 
\begin{equation}   
  g'(z)=\frac12\biggl(\frac1{z}+\frac{\gamma tz^{-2}}{1-tz^{-1}}
-\frac{\gamma t}{1-tz}\biggr)
-\frac{(\gamma-1)(z+1)}{2z(z-t)(z-t^{-1})}\beta(z).
\end{equation} 
Thus we have for $|z|>1$, $z\notin(-\infty,-1)\cup[t^{-1},\infty)$,
\begin{equation}\label{eq-rhp15}
  g(z)= \frac{\gamma}2\log(1-tz)(1-tz^{-1})+\frac12\log z
-\int_{1_{+0}}^z \frac{(\gamma-1)(s+1)}{2s(s-t)(s-t^{-1})} \beta(s)ds
+g_-(1)-\gamma\log(1-t),
\end{equation}
and for $|z|<1$, $z\notin(-1,t]$,
\begin{equation}\label{eq-rhp16}
  g(z)= \frac{\gamma}2\log(1-tz)(1-tz^{-1})+\frac12\log z
-\int_{1_{-0}}^z \frac{(\gamma-1)(s+1)}{2s(s-t)(s-t^{-1})} \beta(s)ds
+g_+(1)-\gamma\log(1-t).
\end{equation}
Now we compute $g_+(1)+g_-(1)$.
From \eqref{eq-g1} with $z=1$ and the evenness of $\psi$, we have 
\begin{equation}   
  g_+(1)+g_-(1)-\pi i 
= \frac{8(\gamma-1)}{\pi} \int_0^{\theta_c} 
\log\biggl(2\sin\frac{\theta}2\biggr) 
\frac{\cos\frac\theta2}{\frac{(1-t)^2}{t}+4\sin^2\frac\theta2} 
\sqrt{\sin^2\frac{\theta_c}2-\sin^2\frac\theta2} d\theta.
\end{equation} 
Substituting $x=(\sin\frac{\theta_c}2)^{-1}\sin\frac\theta2$, 
the above becomes
\begin{equation}\label{eq-rhp18}
\frac{4(\gamma-1)}{\pi}\int_0^1 \log\biggl(2\sin\frac{\theta_c}2x\biggr)
\frac{\sqrt{1-x^2}}{x^2+p^2}dx, 
\qquad p:=\frac{1-t}{2\sqrt{t}\sin\frac{\theta_c}2}
=\frac{\gamma-1}{\sqrt{2\gamma-1}}.
\end{equation}
Using $\frac{2(\gamma-1)}{\pi}\int_0^1 \frac{\sqrt{1-x^2}}{x^2+p^2}dx=1$ 
which is a consequence of the fact that $\psi(\theta)d\theta/(2\pi)$ 
is a probability measure, we have 
\begin{equation}\label{eq-rhp19}   
   g_+(1)+g_-(1)-\pi i
=2\log\biggl(2\sin\frac{\theta_c}2\biggr)
+ \frac{4(\gamma-1)}{\pi}\int_0^1 \log(x)
\frac{\sqrt{1-x^2}}{x^2+p^2}dx.  
\end{equation} 
Now we use $\int_0^1\log x\frac{dx}{\sqrt{1-x^2}}=
\int_0^{\pi/2}\log(\sin\theta)d\theta=-\frac{\pi}2\log2$ to 
rewrite the last term in the above as 
\begin{equation}\label{eq-rhp20}  
  2(\gamma-1)\log2+\frac{4(\gamma-1)(1+p^2)}{\pi}
\int_0^1\frac{\log x}{(x^2+p^2)\sqrt{1-x^2}} dx.
\end{equation} 
Also one can show by residue calculations that 
\begin{equation}   
  \int_0^1\frac{\log x}{(x^2+p^2)\sqrt{1-x^2}} dx
=-\frac{\pi}2\int_1^\infty \frac{dx}{(x^2+p^2)\sqrt{x^2-1}}
+\frac{\pi\log p}{2p\sqrt{1+p^2}}.
\end{equation} 
But one can directly verify that 
\begin{equation}   
 \int \frac{dx}{(x^2+p^2)\sqrt{x^2-1}}
= \frac1{2p\sqrt{1+p^2}} \log\biggl| 
\frac{p\sqrt{x^2-1}+\sqrt{1+p^2}x}{p\sqrt{x^2-1}-\sqrt{1+p^2}x} 
\biggr| +C.
\end{equation} 
Thus from \eqref{eq-rhp19}, \eqref{eq-rhp20} and 
the definition of $p$ in \eqref{eq-rhp18}, we have 
\begin{equation}   
   g_+(1)+g_-(1) -\pi i
= 2\gamma\log\frac{2(\gamma-1)}{2\gamma-1}
+ \log\frac{(2\gamma-1)(1-t)^2}{4t(\gamma-1)^2},
\end{equation} 
which is equal to $2\gamma\log(1-t)-l$.
Now for $z\in \overline{J}$, from \eqref{eq-rhp15} and \eqref{eq-rhp16}, 
we have 
\begin{equation}\label{eq-rhp24}   
  g_+(z)+g_-(z)= V(z)+\log z -l+\pi i.
\end{equation} 
Thus \eqref{eq-g1} yields that the first variational condition 
\eqref{eq-varcon} is satisfied.
On the other hand, for $z\in\Sigma\setminus\overline{J}$, 
$g_+(z)+g_-(z)$ is equal to the right hand side of \eqref{eq-rhp24} 
plus 
\begin{equation}\label{eq-rhp25}   
  -\int_{C_1} \frac{(\gamma-1)(s+1)}{2s(s-t)(s-t^{-1})} \beta(s)ds, 
\end{equation} 
where $C_1=\{ e^{i\theta} : \theta_c\le\theta\le \arg z\}$ 
oriented from $\xi$ to $z$ if $\arg z>0$, 
and $C_1=\{ e^{i\theta} : \arg z \le\theta\le -\theta_c \}$ 
oriented from $\xi^{-1}$ to $z$ if $\arg z<0$. 
For $\theta_c<\arg z<\pi$, \eqref{eq-rhp25} is equal to 
\begin{equation}   
  -\int_{\theta_c}^{\arg z} 
\frac{\sqrt{2}(\gamma-1)\cos\frac\theta2}{t+t^{-1}-2\cos\theta}
\sqrt{\cos\theta_c-\cos\theta} d\theta <0, 
\end{equation} 
and for $-\pi<\arg z<-\theta_c$, it is equal to 
\begin{equation}   
  -\int_{\arg z}^{-\theta_c} 
\frac{\sqrt{2}(\gamma-1)\cos\frac\theta2}{t+t^{-1}-2\cos\theta}
\sqrt{\cos\theta_c-\cos\theta} d\theta <0.
\end{equation} 
Thus the second variational condition of \eqref{eq-varcon} 
is satisfied and the inequality is strict.
\end{itemize}
\end{proof}

For fixed $0<t<1$, we define a weight function on the unit circle $\Sigma$ by 
\begin{equation}   
  \varphi(z;N):= (1-tz)^{-N}(1-tz^{-1})^{-N}.
\end{equation} 
Let $\pi_n(z;N)=z^n+\cdots$ be the $n^{\text{th}}$ monic orthogonal polynomial 
with respect to the measure $\varphi(z;N)dz/(2\pi iz)$ 
on the unit circle, and let $N_n(N)$ be the norm of $\pi_n(z;N)$ : 
\begin{equation}   
  \int_\Sigma \pi_n(z;N)\overline{\pi_m(z;N)} 
(1-tz)^{-N}(1-tz^{-1})^{-N}\frac{dz}{2\pi iz} 
=N_n(N)\delta_{nm}.
\end{equation} 
We also define 
\begin{equation}   
  \pi^*_n(z;N)=z^n\pi_n(z^{-1};N).
\end{equation} 

Define a $2\times 2$ matrix-valued function $Y$ of $z\in\C\setminus\Sigma$ 
by
\begin{equation}\label{as4}
    Y(z;n;N):=\begin{pmatrix}
    \pi_{n}(z;N)&
\int_{\Sigma}\frac{\pi_{n}(s;N)}{s-z}\frac{\varphi(s;N)ds}{2\pi i s^{n}}\\
    -N_{n-1}(N)^{-1}\pi^{*}_{n-1}(z;N)& -N_{n-1}(N)^{-1}
\int_{\Sigma}\frac{\pi^{*}_{n-1}(s;N)}{s-z}\frac{\psi(s;N)ds}{2\pi i s^{n}}
       \end{pmatrix}, \qquad n\geq 1.
\end{equation}
Then $Y(\cdot\thinspace;n;N)$ solves the following
Riemann-Hilbert problem (RHP) (see Lemma 4.1 in \cite{BDJ}) :
\begin{equation}\label{as-RHP}
 \begin{cases}
     Y(z;n;N) \quad\text{is analytic in}\quad z\in\C\setminus\Sigma,\\
     Y_+(z;n;N)=Y_-(z;n;N) \begin{pmatrix}
1&\frac{1}{z^{n}}\varphi(z;N)\\0&1 \end{pmatrix},
\quad\text{on}\quad z\in\Sigma,\\
     Y(z;n;N) \biggl(\begin{smallmatrix} z^{-n}&0\\0&z^n
\end{smallmatrix} \biggr)
=I+O(\frac1{z}) \quad\text{as}\quad z\to\infty.
 \end{cases}
\end{equation}
Here the notation $Y_+(z;n:N)$ (resp., $Y_-$) denotes the limit 
of $Y(z';n;N)$ as $z'\to z$ satisfying $|z'|<1$ (resp., $|z'|>1$).
Note that $n$ and $N$ play the role of external parameters in the above RHP.
One can easily show that the solution of the above RHP is unique,
hence \eqref{as4} is the unique solution of the above RHP.
This RHP formulation of orthogonal polynomials on the unit circle
is an adaptation of a result of Fokas, Its and Kitaev in \cite{FIK}
where orthogonal polynomials on the real line are considered.

From \eqref{as4}, it is easy to check that 
\begin{eqnarray}
\label{as7.7}   N_{n-1}(N)^{-1} &=& -Y_{21}(0;n;N), \\
\label{as7.8}  \pi_n(z;N) &=& Y_{11}(z;n;N), \\
\label{as7.9}  \pi^{*}_{n}(z;N) &=& z^nY_{11}(z^{-1};n;N)
= Y_{21}(z;n+1;N)(Y_{21}(0;n+1;N))^{-1}.
\end{eqnarray}
The asymptotics of the above quantities can be obtained by 
applying Deift-Zhou method for Riemann-Hilbert problem \eqref{as-RHP}.
A reference for Deift-Zhou method is \cite{deift}.
In \cite{BDJ} and \cite{BR2}, similar asymptotics are obtained 
for different weight function $e^{t(z+z^{-1})}$ as $t,n\to\infty$.
It is interesting to compare the following results with 
Proposition 5.1 in \cite{BR2}.

\begin{prop}\label{prop-rhp}
For each $n$ and $N$, define $x\in\R$ by 
\begin{equation}\label{eq-as37}   
  \frac{2t}{1+t}\frac{N}{n}=
1-\biggl[\frac{1-t}{2(1+t)^2}\biggr]^{1/3}\frac{x}{n^{2/3}}.
\end{equation} 
Also let 
\begin{equation}   
  l:= \frac{2N}{n}\log\frac{(2N-n)(1-t)}{2(N-n)}
-\log\frac{n(2N-n)(1-t)^2}{4t(N-n)^2}
\end{equation}
and for $\frac{N}{n}>\frac{1+t}{2t}$, let 
\begin{equation}   
  \sin\frac{\theta_c}2:= \frac{(1-t)}{2(N-n)}\sqrt{\frac{n(2N-n)}{t}}.
\end{equation} 
There exits $M_0>0$ such that as $n,N\to\infty$, the following asymptotic 
results hold.
\begin{enumerate}
\item If $0\le \frac{2t}{1+t}\frac{N}{n} \le a$ for $0<a<1$,
then 
\begin{equation}   
   |N_{n-1}(N)^{-1}-1|,\quad  |\pi_n(0;N)| \le Ce^{-cn},
\end{equation} 
for some constants $C$, $c$ which may depend on $a$.
\item If $a\le\frac{2t}{1+t}\frac{N}{n} \le 1- Mn^{-2/3}$ 
for some $M>M_0$ and $0<a<1$, then 
\begin{equation}   
   |N_{n-1}(N)^{-1}-1|,\quad  |\pi_n(0;N)| \le 
\frac{C}{n^{1/3}} e^{-cx^{3/2}}
\end{equation} 
for some constant $C$ and $c$ depending on $M$.   
\item If $-M\le x\le M$ for some $M>0$, then 
\begin{equation}   
   \biggl|N_{n-1}(N)^{-1}-1
-\biggl[\frac{2(1+t)^2}{1-t}\biggr]^{1/3}\frac{v(x)}{n^{1/3}}\biggr|,
\quad  \biggl|\pi_n(0;N)
+(-1)^n\biggl[\frac{2(1+t)^2}{1-t}\biggr]^{1/3}\frac{u(x)}{n^{1/3}}\biggr| 
\le \frac{C}{n^{2/3}}
\end{equation} 
for some constant $C$ depending on $M$.
\item If $1+Mn^{-2/3}\le \frac{2t}{1+t} \le a$ for some 
$M>M_0$ and $a>1$, then 
\begin{equation}   
  \biggl| \frac{e^{-nl}}{\sin\frac{\theta_c}2}N_{n-1}(N)^{-1}-1\biggr|,
\quad \biggl| \frac{(-1)^n}{\cos\frac{\theta_c}2}\pi_n(0;n;N)-1\biggr|
\le \frac{C}{\frac{2t}{1+t}N-n}
\end{equation}
for some constant $C$ depending on $M$.
\item If $a\le \frac{2t}{1+t}$ for some $a>1$, 
\begin{equation}   
  \biggl| \frac{e^{-nl}}{\sin\frac{\theta_c}2}N_{n-1}(N)^{-1}-1\biggr|,
\quad \biggl| \frac{(-1)^n}{\cos\frac{\theta_c}2}\pi_n(0;n;N)-1\biggr|
\le \frac{C}{n}
\end{equation} 
for some constant $C$ depending on $a$.
\end{enumerate}
\end{prop}

\begin{rem}
From the calculations analogous to Section 10 of \cite{BR2}, 
in addition to the above asymptotics results, 
we can obtain more results similar to 
those in Section 5 of \cite{BR2}.
For example, suppose $x$ defined in \eqref{eq-as37} above satisfies 
$c_1\le x\le c_2$ for some constants $c_1,c_2\in\R$ 
(hence we are in the case (iii) 
of above proposition). 
For $\alpha>0$, define $w\in\R$ by 
\begin{equation}
  \alpha=1-\biggl[\frac{2(1+t)^2}{1-t}\biggr]^{1/3}\frac{2w}{n^{1/3}}.
\end{equation}
Then we have for $w>0$ fixed, 
\begin{eqnarray}
  \lim_{n\to\infty} (-1)^n(1+t\alpha)^{-N} \pi_n(-\alpha;N)
&=&-m_{12}(-iw;x),\\
  \lim_{n\to\infty} (1+t\alpha)^{-N} \pi_n^*(-\alpha;N)
&=&m_{22}(-iw;x),
\end{eqnarray}
and for $w<0$ fixed, 
\begin{eqnarray}
  \lim_{n\to\infty} (-1)^n(1+t\alpha)^{-N} \pi_n(-\alpha;N)
&=&m_{11}(-iw;x)e^{\frac83w^3-2xw},\\
  \lim_{n\to\infty} (1+t\alpha)^{-N} \pi_n^*(-\alpha;N)
&=&-m_{21}(-iw;x)e^{\frac83w^3-2xw},
\end{eqnarray}
where $m(z;x)$ is the solution to the Riemann-Hilbert problem 
for Painlev\'e II equation with special monodromy data 
$p=-q=1$, $r=0$ (see, e.g. (2.15) of \cite{BR2}).
These results are parallel to (5.21), (5.22), (5.25), (5.26) 
of \cite{BR2}.
But in this paper, we only need Proposition \ref{prop-rhp} above.
\end{rem}

We are not going to present the detail of the proof because the computation 
is parallel to that of Lemma 5.1 and Lemma 6.3 of \cite{BDJ} 
(see also Proposition 5.1 of \cite{BR2}) where the authors consider 
the same asymptotic problem with 
different weight function $e^{\sqrt{\lambda}\cos\theta}$.
Instead we give some indication why we have Painlev\'e II function 
in the case (iii).
Let us consider only when $x>0$.
Define a $2\times 2$ matrix valued function $m(z)$ by 
\begin{eqnarray}
  m(z):= \begin{cases}
Y(z;n;N)\begin{pmatrix} 0&-(1-tz)^{-N}\\(1-tz)^N&0 
\end{pmatrix}, &|z|<1\\
Y(z;n;N)\begin{pmatrix} z^{-n}(1-tz^{-1})^{-N}&0\\
0&z^n(1-tz^{-1})^N
\end{pmatrix}, &|z|>1,
\end{cases}
\end{eqnarray}
where $Y$ is defined in \eqref{as4}.
From the RHP for $Y$, $m$ solves a new RHP : 
$m$ is analytic in $z\in\C\setminus\Sigma$, 
$m(z) =I+O(\frac1{z})$ as $z\to\infty$, and 
satisfies the jump condition $m_+(z)=m_-(z)v(z)$ on $z\in\Sigma$ where 
\begin{equation}\label{as-newRHP}
     v(z)= \begin{pmatrix}
1&-z^{n}(1-tz)^{-N}(1-tz^{-1})^N\\
z^{-n}(1-tz)^N(1-tz^{-1})^{-N}&0    \end{pmatrix}.
\end{equation}
The choice of $m$ above is related to the equilibrium measure 
in Lemma \ref{as-lem1}. 
The role of equilibrium measure in RHP for orthogonal polynomials 
is discussed in \cite{DKMVZ2, DKMVZ3} (see also \cite{BDJ}).
The two RHP's for $Y$ and $m$ are algebraically related and are equivalent 
in the sense that a solution to one RHP implies a solution to the other RHP.
Note that the jump matrix $v(z)$ has the factorization 
$v(z)=v_-(z)v_+(z)$ where 
\begin{equation}
   v_-(z)= \begin{pmatrix}
1&0\\
z^{-n}(1-tz)^N(1-tz^{-1})^{-N}&1    \end{pmatrix}, 
\qquad 
     v_+(z)= \begin{pmatrix}
1&-z^{n}(1-tz)^{-N}(1-tz^{-1})^N\\
0&1    \end{pmatrix}.
\end{equation}
Hence by usual deformation technique of RHP, we can bring 
the matrix $v_+$ to a contour in $|z|<1$, and the matrix $v_-$ to 
a contour in $|z|>1$.
By the assumption of $N$ and $n$ in case (iii), 
except in a neighborhood of $z=-1$,
we can find a new contour where 
the off diagonal entries of $v_\pm$ decay exponentially as $n\to\infty$. 
Hence the main contribution to the RHP as $n\to\infty$ comes only from 
a neighborhood of $z=-1$.
This is exactly related to the fact that the support of the equilibrium 
measure is Lemma \ref{as-lem1} has special point $z=-1$ 
at which a new gap opens up when 
$\gamma (=\frac{N}{n}) = \frac{1+t}{2t}$. 
Now let us focus on the neighborhood of $z=-1$. 
The (12) entry of $v$ is $-e^{h(z)}$ where 
\begin{equation}\label{eq-as48}
   h(z) = n\log z -N \log(1-tz) +N\log(1-tz^{-1}).
\end{equation}
Set $z=-1+s$. 
Using the definition of $x$ in \eqref{eq-as37}, 
expansion of \eqref{eq-as48} becomes 
\begin{equation}
\begin{split}
  h(z)&= n\log(-1) - x\biggl[\frac{1-t}{2(1+t)^2}\biggr]^{1/3}(n^{1/3}s)
- \frac{x}{2n^{1/3}} \biggl[\frac{1-t}{2(1+t)^2}\biggr]^{1/3}(n^{1/3}s)^2
+ \frac{1-t}{6(1+t)^2} (n^{1/3}s)^3 + \cdots\\
&\sim n\log(-1) - 2x\tilde{s}+\frac83\tilde{s}^3 
\qquad \text{as $n\to\infty$,}
\end{split}
\end{equation}
where $\tilde{s}:=\frac12\bigl[\frac{1-t}{2(1+t)^2}\bigr]^{1/3}(n^{1/3}s)$.
Thus we are lead to a RHP with the jump matrix 
\begin{equation}
     \tilde{v}(z)= \begin{pmatrix}
1& -(-1)^ne^{2(-x\tilde{s}+\frac43\tilde{s}^3)}\\
(-1)^ne^{-2(-x\tilde{s}+\frac43\tilde{s}^3)}&0    \end{pmatrix}
\end{equation}
on the line $i\R$ oriented from $+i\infty$ to $-i\infty$.
After rotation by $-\pi/2$, this is precisely the jump matrix 
for the Painlev\'e II equation with parameter $p=-q=1$, $r=0$ 
(see, e.g. (2.15) of \cite{BR2} : the term $(-1)^n$ in off diagonal 
entries can be simply conjugated out).
Thus the $m$, therefore $Y$, can be expressed in terms of Painlev\'e II 
solution in the limit $n\to\infty$ in the case (iii) with $x>0$.

\bibliographystyle{plain}
\bibliography{paper5}

\end{document}